\documentstyle{amsppt}
\loadbold
\magnification=\magstep1
\hsize=6.5 true in
\vsize=8.9 true in

\overfullrule=0pt
\def\eqnum#1{\eqno (#1)}

\def\Aff{\operatorname{Aff}}
\def\Aut{\operatorname{Aut}}
\def\const{\operatorname{const}}

\def\hull{\operatorname{hull}}
\def\id{\operatorname{id}}
\def\Im{\operatorname{Im}}

\def\Ker{\operatorname{Ker}}

\def\pr{\operatorname{pr}}
\def\Re{\operatorname{Re}}
\def\rk{\operatorname{rk}}
\def\supp{\operatorname{supp}}
\def\SL{\operatorname{SL}}
\def\normo#1{\left\| #1 \right\|}
\def\modo#1{\left| #1 \right|}
\def\eqnum#1{\eqno (#1)}

\font\bfit=cmbxti10
\font\smallit=cmti8
\topmatter
\title{Liouville and Carath\'eodory coverings in Riemannian
and complex geometry}\endtitle

\author{Vladimir Ya. Lin, Mikhail Zaidenberg}\endauthor

\thanks{The authors thank the Max-Planck-Institut f\"ur Mathematik in Bonn
for hospitality. The research of the first author was partially
supported by the Fund of the Israel Science Foundation.}
\endthanks

\address{Department of Mathematics, Technion, Haifa 32000, Israel}\endaddress

\email{vlin\@techunix.technion.ac.il}\endemail

\address{Universit\'{e} Grenoble I,
Institut Fourier et Laboratoire de Math\'ematiques
associ\'e au CNRS,
BP 74,
38402 St. Martin d'H\`{e}res--c\'edex,
France}\endaddress

\email{zaidenbe\@puccini.ujf-grenoble.fr}\endemail

\leftheadtext{Vladimir Lin, Mikhail Zaidenberg}
\rightheadtext{Liouville and Carath\'eodory coverings}
\endtopmatter

\document
\hfill{\smallit To S. G. Krein with admiration and love.}

\head{Introduction}\endhead

A Riemannian manifold resp. a complex space $X$
is called {\it Liouville} if it carries no nonconstant
bounded harmonic resp. holomorphic functions.
It is called {\it Carath\'eodory}, or {\it Carath\'eodory  hyperbolic},
if bounded harmonic resp. holomorphic functions separate the points of $X$.
The problems which we discuss in this paper arise from the following
question:
\vskip0.12pc

{\sl When a Galois covering $X$ with Galois group $G$ over
a Liouville base $Y$ is Liouville or, at least, is not Carath\'eodory
hyperbolic?}
\vskip0.12pc

An infinite abelian\footnote
{A Galois covering $X\to Y$ with Galois group $G$
is referred to as a $G$-covering. If $G$ is abelian
(resp. nilpotent, solvable, polycyclic, etc.)
the covering is called {\it abelian}
(resp. {\it nilpotent, solvable, polycyclic}, etc.).
It is said to be {\it cocompact} if its base $Y$ is compact.}
covering of a Liouville base $Y$ need not be Liouville even
for an open Riemann surface $Y$.
In \cite{LySu} a $\Bbb Z^{\infty}$-covering of this kind
was constructed.
Moreover, there is a non-Liouville $\Bbb Z$-covering
of a Liouville complex surface
\cite{Li2} (see Remark 1.9.1).
Thus, to ensure the Liouville property
of $X$ one must subject $Y$ to a stronger condition.

By this reason, we require $Y$ to be compact or, more generally, to
carry no nonconstant bounded subharmonic resp.
plurisubharmonic functions. Then, to some extent, the coverings over
Riemannian and complex spaces behave similarly.
Roughly speaking, $X$ is Liouville if $G$
is small enough, say nilpotent \cite{LySu,\, Li2};
and a solvable cocompact covering
can be even Carath\'eodory hyperbolic (see
Theorem 1.6 and \S 3).
But in the intermediate class of polycyclic coverings
this similarity fails: such a covering
over a compact Riemannian resp. K\"ahler base $Y$
is Liouville \cite{Ka1},
while there is a non-K\"ahler compact complex surface
with non-Liouville polycyclic universal covering \cite{Li2}
(see Theorem 1.1 and \S 4).

We start in \S 1 with a brief survey of some known results,
sketching a few proofs, and proceed with certain new observations.
In particular, combining a theorem of Varopoulos \cite{VSCC} with
a theorem in \cite{LySu},
we establish that a $G$-covering over a compact Riemannian resp. K\"ahler
manifold is Liouville if $G$ is an
extension of an almost nilpotent group by $\Bbb Z$ or by $\Bbb Z^2$
(see Theorems 1.4, 1.6, and Corollary 1.8).

We discuss Liouville-type properties not only for
$G$-coverings $X\to Y$ but also for general $G$-spaces
$X$, Riemannian or complex.
We show, in particular, that $X$ must be Liouville whenever
the induced diagonal $G$-action in $X\times X$
has a dense orbit (Proposition 1.10$(b)$).

Further, we consider the period subgroup of
$G$ consisting of all $g\in G$ that do not affect
the bounded harmonic resp. holomorphic functions on $X$.
If the given $G$-action $T$ on $X$ is cocompact,
or $T$ is {\it ultra-Liouville} (see \S 1.5 below for the definition)
and $G$ is amenable, then the period subgroup contains any
central element of $G$ and,
moreover, any element with finite
conjugacy class \cite{Li2} (see Corollary 1.14).
In \S 2 we extend the latter result to the
FC-hypercentral elements and establish Liouville property of
FC-(hyper)nilpotent coverings (see Definitions 2.1 for terminology).
\vskip0.14pc

\S 3 and \S 4 are devoted to certain
examples of non-Liouville and,
especially, Carath\'eodory hyperbolic cocompact coverings
with relatively small Galois groups. In \S 3 for any compact
Riemann surface $Y$ of genus $g \ge 2$
we produce a Carath\'eodory hyperbolic
two-step solvable covering $X$ over $Y$.
This is based on a construction due to Lyons and Sullivan \cite{LySu}.
\vskip0.14pc

In \S 4 we consider the universal
covering $X\to \Cal I$ over an Inoue surface $\Cal I$ \cite{In},
which is a non-K\"ahler compact complex surface with
a polycyclic fundamental group $G=\pi_1(\Cal I)$. In fact,
$X\cong \Bbb H\times \Bbb C$ \ ($\Bbb H$ is the upper halfplane);
it is neither Liouville nor Carath\'eodory.
We show that $X$ admits bounded holomorphic functions which are
nonconstant on the orbits of suitable infinite conjugacy classes
in $G$. Notice that, by Corollary 1.14,
this is impossible for a finite conjugacy class.
\smallskip

Throughout the paper, all manifolds and complex spaces are assumed to be
connected.

\head{\S 1. Liouville-type properties of coverings and $G$-spaces:
A survey}\endhead

\noindent This brief survey is neither complete nor chronological;
it contains only some selected results on Liouville-type properties.
We do not address the case of harmonic functions on
a discrete group with a probability measure, which is closely
related to our topic
(see e. g. \cite{Av, Fu, Mar, KaVe, Ka1, Ka3, VSCC, Wo}).
\bigskip

\noindent{\bfit Coverings over a compact base}
\smallskip

\proclaim{1.1. Theorem {\rm \cite{LySu,\, Ka1}}} Let $X \to Y$ be
a Galois covering with Galois group $G$ over
a Riemannian manifold $Y$.
Suppose that the base \ $Y$ is compact; then
\vskip0.14pc

\noindent $a)$ $X$ is Liouville whenever $G$ is polycyclic\footnote
{i. e. $G$ admits a finite normal series
with cyclic quotients, or, equivalently,
$G$ is solvable and all its subgroups are finitely generated
(see e. g. \cite{Ha,\, Se}; in \cite{Ha} polycyclic groups were called
{\it supersolvable}).} or of subexponential growth\footnote
{Any finitely generated solvable nonpolycyclic group is of
exponential growth (see [Mi]).} \cite{Ka1};
\vskip0.14pc

\noindent $b)$ if $G$ is nilpotent, then $X$ carries no nonconstant positive
harmonic functions\footnote{See also \cite{Gui,\, Mar,\, LiPi}.}
\cite{LySu}.
\smallskip

\noindent Furthermore, for any base $Y$ the group $G$ must be amenable\footnote
{According to von Neumann, $G$ is called {\it amenable}
if the Banach space $L^\infty(G)$ of all bounded functions on $G$
admits a $G$-(right)invariant mean; see, for instance, \cite{Gre}.}
if $X$ is Liouville \cite{LySu}.
\endproclaim

\remark{\bf 1.2. Remarks} 1. Every compact Riemann
surface $R$ of genus $g\ge 2$ admits a non-Liouville solvable
covering \cite{LySu} (see also \S 3 for a stronger example of such kind).
Thus, Theorem 1.1$(a)$, in general,
does not hold for nonpolycyclic solvable cocompact
coverings.
\vskip0.14pc

\noindent 2. Let $X \overset{G}\to{\longrightarrow} Y
{\overset{K}\to{\longrightarrow}} Z$
be the tower of Galois coverings corresponding to a
group extension $1\to G\to \widetilde G\to K\to 1$.
If $K$ is finite then $Y=X/G$ is compact.
Hence, the statements $(a)$, $(b)$ of Theorem 1.1 hold true for any finite
extension $\widetilde G$ of $G$ if $G$ is as in these statements.

\noindent 3. One says that a group $G$ is {\it almost nilpotent}\footnote
{or {\it virtually nilpotent},
or also {\it nilpotent--by--finite}.} (resp.
{\it almost solvable, almost polycyclic}, etc.)
if it contains a nilpotent (resp. solvable, polycyclic,
etc.) subgroup of finite index. Such a subgroup
may clearly be assumed being normal. Thus, by Theorem 1.1$(a),(b)$,
an almost polycyclic resp. an almost
nilpotent covering over a compact Riemannian manifold is Liouville resp.
carries no nonconstant positive harmonic function.
\vskip0.13pc

\noindent 4. A holomorphic function on a K\"ahler manifold is harmonic
with respect to the Laplace-Beltrami operator
related to the K\"ahler metric. Hence, Theorem 1.1$(a)$ holds true for
holomorphic functions on coverings of compact K\"ahler manifolds. It holds
also true for compact {\it semi-K\"ahlerian} manifolds \cite{Ka2}.
The latter class of Hermitian manifolds may actually be characterized by the
harmonicity of all holomorphic germs \cite{Ka2} (cf. \cite{Ga}).
\vskip0.13pc

\noindent 5. Theorem 1.1$(a)$ does not hold, in general,
in the case of holomorphic functions on polycyclic coverings of
non-K\"ahler compact complex manifolds.
For instance, consider the universal covering $X$ of the Inoue surface
$\Cal I$ (\cite{In, Li2}; see also \S 4).
The Galois group $G=\pi_1(\Cal I)$ contains a normal subgroup
$G_0\cong\Bbb Z^3$ with $G/G_0\cong\Bbb Z$.
Hence, $G$ is a semidirect product $\Bbb Z^3\leftthreetimes \Bbb Z$
and so it is a metabelian (i. e. a two-step solvable) polycyclic group.
Furtermore, $X\cong\Bbb H\times\Bbb C$ \ ($\Bbb H$ is the upper halfplane)
is not Liouville. However, nilpotent coverings of compact complex spaces
are Liouville (\cite{Li2}; see Theorem 1.6 below).
\vskip0.13pc

\noindent 6. The last assertion of Theorem 1.1 has no direct analog in
complex geometry. For instance, set
$U=\{(x,\,y)\in\Bbb C^2\mid \,\text{either} \ \modo{x}<1 \
\text{and} \ y\ne 2 \ \text{or} \ x\ne 2 \ \text{and} \ \modo{y}<1\}$.
Then $\pi_1(U)$ is a free group of rank $2$ (i. e. nonamenable),
and the universal covering of $U$ is Liouville. Moreover,
for any finitely presented group $G$ there is
a Stein manifold $Y$ with $\pi_1(Y)\cong G$ and Liouville
universal covering $X$ (V. Lin, unpublished; cf. also \cite{Li1, \S 7.2}).
\endremark

\definition{1.3. Definition} A Riemannian manifold $Y$ is called
{\it transient} if it carries a nonconstant bounded
subharmonic function, or, equivalently,
if it possesses positive Green function.
Nontransient manifolds are called {\it reccurent};
this property is equivalent to the reccurence
of the random motion on $Y$
(see e. g. \cite{SNWC, Gri, LySu}).
In \cite{Li2} the following terminology was suggested: a Riemannian manifold
resp. a complex space $Y$ is called {\it ultra-Liouville}
if any bounded {\sl continuous} subharmonic resp. plurisubharmonic function
on $Y$ is constant. Ultra-Liouville Riemannian
manifolds are recurrent, and vice versa.

Any connected Zariski open subset $Y$ of a compact complex space
$\overline Y$
(for instance, any quasiprojective complex variety $Y$) is ultra-Liouville.
Indeed, by a theorem of Grauert and Remmert \cite{GraRe}
(see also \cite{BoNa})
every bounded plurisubharmonic function on $Y$ admits
a plurisubharmonic extension to $\overline Y$
and, hence, by the maximum principle, it is constant. Note that a smooth
quasiprojective complex variety, being endowed with a Riemannian metric,
may be transient; e. g., this is so for
$Y = \Bbb C^n$, \ $n \ge 2$, with its Euclidean metric.
\enddefinition

The following recurrence criterion of cocompact coverings
was proved in \cite{VSCC, X.3}\footnote{
See the references therein. For abelian coverings
this theorem was proved in \cite{Gui} and \cite{LySu}. For the classical
case of Riemann surfaces
see e. g. \cite{My, Ne, Roy, Mo, Ts}.}.

\proclaim{1.4. Theorem} Let $X\to Y$
be a Galois covering with Galois group $G$
over a compact Riemannian manifold $Y$.
Then $X$ is recurrent $($or, which is equivalent, ultra-Liouville$)$
if and only if $G$ is a {\it Varopoulos group}, that is,
a finite extension of
one of the groups $\bold 1$, $\Bbb Z$, and $\Bbb Z^2$.
\endproclaim

\noindent{\bfit Ultra-Liouville actions and
coverings over a noncompact base}
\smallskip

\definition{1.5. Notation and definitions}
Given a Riemannian manifold resp. a complex space $X$,
we denote by $I(X)$ the group of all its homotheties\footnote
{By a homothety of
a Riemannian manifold $(X, \,d)$ we mean
a transformation $g\colon \,X \to X$ such that $d(gx,\,gy) \equiv C d(x,\,y)$
with some constant $C = C(g)$ which
does not depend on $x,\,y \in X$.}
${\text {Homo}}\,(X)$ resp.
the group of all its biholomorphic automorphisms $\text{Aut}\,(X)$. By
${\Cal H}={\Cal H}(X)$ we denote
the space $Harm^{\infty}(X)$ resp. $H^{\infty}(X)$
of all bounded complex valued harmonic resp. holomorphic functions
on $X$.

Clearly, $I(X)$ acts in ${\Cal H}(X)$.
We say that the action of a subgroup $G\subseteq I(X)$
on $X$ is {\it ultra-Liouville}
if $X$ admits no nonconstant $G$-invariant bounded continuous
subharmonic resp. plurisubharmonic functions.
If the quotient $Y = X/G$ exists in the same category, then the
$G$-action on $X$ is ultra-Liouville if and only if $Y$ is ultra-Liouville
in the sense of Definition 1.3.
\smallskip

Let $Z(G)$ denote the center of a group $G$.
Consider the upper central series of $G$
$$
{\bold  1}=Z_0(G) \vartriangleleft Z(G) = Z_1(G)  \vartriangleleft
Z_2(G) \vartriangleleft \dots \vartriangleleft Z_n(G) \vartriangleleft
\dots \vartriangleleft G\,;
$$
here $Z_n(G)$ is the total preimage $p_{n-1}^{-1}(Z(G/Z_{n-1}(G)))$
of $Z(G/Z_{n-1}(G))$
under the natural surjection $p_{n-1}\colon \,G\to G/Z_{n-1}(G)$,
\ $n=1,2,\dots$ \,.
The upper central series is continued transfinitely in the usual way,
by defining $Z_\alpha(G)=\bigcup_{\beta<\alpha} Z_\beta(G)$, when
$\alpha$ is a limit ordinal.

The group $G$ is called {\it $\omega$-nilpotent} if it coincides with
the union $Z_{\omega}(G) = \bigcup_{n \in \Bbb N} Z_n(G)$.
$\,\,G$ is called {\it hypernilpotent}, or also {\it hypercentral}, if
$G= Z_{\lim}(G)$, where $Z_{\lim}(G) = \bigcup_{\alpha} Z_\alpha (G)$
is the {\it hypercenter} of $G$ (here $\alpha$ runs over all the ordinals).
\enddefinition

The following theorem was proved for $\omega$-nilpotent coverings
of Riemannian manifolds in \cite{LySu}, and in its present
form in \cite{Li2}, by different methods.

\proclaim{1.6. Theorem {\rm \cite{LySu, Li2}}} Let $X$ be a Riemannian
manifold resp. a complex space, and let $G$ be a hypernilpotent subgroup
of $I(X)$. The space $X$ is Liouville whenever
the $G$-action on $X$ is ultra-Liouville. In particular, if $X\to Y$ is
a hypernilpotent covering over an ultra-Liouville
$($Riemannian or complex$)$ base $Y$,
then $X$ is Liouville.
\endproclaim

\remark{\bf 1.7. Remark} By the maximum principle, any cocompact
$G$-action\footnote
{That is, a $G$-action such that $GT = X$ for some
compact set $T\subseteq X$.}
on $X$ is ultra-Liouville. Hence, for hypernilpotent
coverings over a compact Riemannian manifold $Y$ the last assertion of
Theorem 1.6 follows from Theorem 1.1$(a)$ (but not vice versa!). Indeed,
being a quotient of a finitely generated group $\pi_1(Y)$, the Galois
group of a Galois covering $X\to Y$ is
finitely generated, too. But a finitely generated
hypernilpotent group is nilpotent and polycyclic
(see \cite{Ha, Se} resp. Remark 2.2.1 and references therein
for the case of finitely generated
$\omega$-nilpotent resp. finitely generated hypernilpotent groups).

However, unlike Theorem 1.1$(a)$, Theorem 1.6 applies to complex
spaces (K\"ahler or not) and also to ultra-Liouville
actions, which may be neither free nor properly discontinuous nor cocompact.
\endremark
\medskip

From Theorems 1.4, 1.6 we obtain such a corollary.

\proclaim{1.8. Corollary} Let $X \to Y$ be a Galois covering
with Galois group $G$ over a compact Riemannian resp.
K\"ahler manifold $Y$. If $G$ is an extension of an
almost hypernilpotent group
by a Varopoulos group\footnote
{see Theorem 1.4.}, then $X$ is Liouville.
\endproclaim

\remark{\bf 1.9. Remarks} 1. Corollary 1.8 does not apply to general
compact complex manifolds. Indeed \cite{Li2}, let $X \to \Cal I$
be the universal covering over the Inoue surface $\Cal I$
(see Remark 1.2.3 and \S 4).
The semidirect decomposition
$G \cong \Bbb Z^3 \leftthreetimes \Bbb Z$ provides the
tower of Galois coverings
$X {\overset{\Bbb Z^3}\to{\longrightarrow}} Y
{\overset{\Bbb Z}\to{\longrightarrow}} \Cal I$.
If $Y$ were ultra-Liouville,
then, by Theorem 1.6, the abelian covering
$X {\overset{\Bbb Z^3}\to{\longrightarrow}} Y$
would be Liouville, which is wrong.
Hence, due to Theorem 1.6, $Y$ is a Liouville, but not ultra-Liouville
$\Bbb Z$-covering over a compact complex surface $\Cal I$.
Furthermore, the above covering
$X {\overset{\Bbb Z^3}\to{\longrightarrow}} Y$ produces a tower
of three $\Bbb Z$-coverings \
$X=X_1 {\overset{\Bbb Z}\to{\longrightarrow}} X_2
{\overset{\Bbb Z}\to{\longrightarrow}} X_3
{\overset{\Bbb Z}\to{\longrightarrow}} X_4=Y$.
Since $X_4=Y$
is Liouville and $X_1=X$ is not Liouville, at least one
of $X_i$, \ $1\le i\le 3$,
is a non-Liouville $\Bbb Z$-covering over a Liouville base.
\smallskip

\noindent 2. Theorem 1.4 does not hold for coverings over
a noncompact ultra-Liouville
Riemannian manifold $Y$.
Consider, for instance, the maximal abelian covering
$X \to Y$ over the punctured Riemann sphere
$Y = \Bbb P^1 \setminus \{3 \ \text{points}\}\cong\Bbb C\setminus\{0\,,1\}$.
The Riemann surface $X$
can be realized as an analytic curve in $\Bbb C^2$, namely, the curve with
the equation $e^x + e^y = 1$. The covering projection
$X \to Y \cong \Bbb C \setminus \{0,\,1\}$ is
$(x,\,y) \longmapsto e^x$. The Galois group $G$ of this covering is
isomorphic to $H_1(Y;\,\Bbb Z) \cong \Bbb Z^2$.
It is known \cite{McKSu, LyMcK} that $X$ is transient; hence, $X$
is not ultra-Liouville whereas $G\cong \Bbb Z^2$ is a Varopoulos group.
Note that $X$ in this example is Liouville
(see \cite{Dem, Wa, Sh} or Theorem 1.6).
\endremark
\medskip

The next proposition contains some new observations
concerning ultra-Liouville actions.

\proclaim{1.10. Proposition} $a)$ Let the action
of a subgroup $G\subseteq I(X)$ on a Riemannian manifold resp.
on a complex space $X$ be ultra-Liouville. Then any $G$-orbit in $X$
is a uniqueness set for the function space $\Cal H=\Cal H(X)$.
I. e., if $h\in {\Cal H}$, \ $x_0 \in X$, and
$h\mid Gx_0 = 0$, then $h = 0$.
\smallskip

\noindent $b)$ If the induced diagonal $G$-action
$g\colon \,(x,\,y) \mapsto (gx,\,gy)$ on $X \times X$
is ultra-Liouville, then $X$ is Liouville.
\endproclaim

The proof will be done in \S 1.18.

\remark{\bf 1.11. Remark} 1. For complex spaces the following
stronger form of $(a)$
was proven in \cite{BoNa}: a bounded holomorphic
function on a complex space $X$ with an ultra-Liouville $G$-action
is constant whenever the set of its values on some $G$-orbit is finite.
\vskip0.14pc

\noindent 2. The complement $Gx_0\setminus K$
of any finite subset $K\subset Gx_0$ also is a uniqueness set for
$\Cal H$.
\vskip0.14pc

\noindent 3. It follows from Proposition 1.10$(b)$ that
a complex space $X$ is Liouville if the action of the group
$I(X)$ on $X$ is almost doubly transitive, meaning
that the induced diagonal $G$-action on $X \times X$ possesses a dense orbit.
This simple observation yields yet another proof of the classical
Liouville Theorem. Indeed, the affine transformation group
$\Aff\,(\Bbb C)= \Aut\,(\Bbb C)$ is doubly transitive on $\Bbb C$; hence,
the diagonal action on $\Bbb C\times\Bbb C$ is transitive outside of
the diagonal $\Delta$, and $(\Bbb C\times\Bbb C)\setminus\Delta$ is
a dense orbit in $\Bbb C\times\Bbb C$.
\endremark

\definition{1.12. Definition} Given a $G$-space $X$,
the corresponding $G$-action in the vector space $\Bbb C^X$
of all complex valued functions on $X$ is denoted by $f\mapsto f^g$, \
$f^g(x)=f(gx)$.
We say that an element $g \in G$ is a {\it period}
of a function $f \in \Bbb C^X$, or $f$-{\it period}, if $f$ is $g$-invariant,
i. e. $f(gx) = f(x)$ for all $x \in X$.
For a function $f\in \Bbb C^X$ the set of all its periods form a subgroup
in $G$, which is denoted by $G_f$. It is a stationary subgroup of $f$
with respect to the $G$-action on $\Bbb C^X$.
For a subspace ${\Cal F} \subseteq \Bbb C^X$ denote by $G_{\Cal F}$
the intersection of all the subgroups $G_f$, \ $f \in \Cal F$.
We call $G_{\Cal F}$ the {\it subgroup of $\Cal F$-periods},
or simply the {\it period subgroup}.
It is easily seen that $G_{\Cal F}$ is a normal subgroup of
$G$ if $\Cal F$ is $G$-invariant. In particular,
for any subgroup $G\subseteq I(X)$ the $\Cal H$-period subgroup
$G_{\Cal H}$ is normal in $G$, where, as before, $\Cal H = \Cal H (X)$
denotes the space of bounded harmonic resp. holomorphic functions on a
Riemannian manifold resp. on a complex space $X$.

For a subgroup $G\subseteq I(X)$ and an element $s\in I(X)$
we denote by $[s,\,G]$ the subgroup of $I(X)$ generated by all the
commutators $[s,\,g] = sgs^{-1}g^{-1}$, \ $g \in G$.
\enddefinition

The following theorem provides certain information on
the period subgroup $I(X)_f$ of any bounded harmonic resp.
holomorphic function $f$ on $X$.

\proclaim{1.13. Theorem {\rm \cite{Li2, Thms. 2.10, 3.9}}}
Let, as before, $X$ be a Riemannian manifold resp.
a complex space, and let $G$ be a subgroup of the group
$I(X)={\text{Homo}} \,(X)$ resp. $I(X)={\text{Aut}} \,(X)$.
Suppose that one of the following two conditions is fulfilled:
\smallskip

\noindent $\,*$ $G$ is amenable and its action on $X$ is ultra-Liouville;
\smallskip

\noindent $\,*$ the action of $G$ on $X$ is cocompact.
\smallskip

\noindent Let $f$ be a bounded harmonic resp. holomorphic function on $X$.
Suppose that an element $s\in I(X)$ satisfies the condition
$[s, G]\subseteq I(X)_f$. Then $s\in I(X)_f$.
In other words, if $f$ is invariant under all the commutators
$[s,\, g]=sgs^{-1}g^{-1}$, \ $g\in G$, then $f$ is also $s$-invariant.
\endproclaim

In Corollary 1.14 and Remarks 1.15 below we use the same notation
as in Theorem 1.13; see also Definitions 1.5, 1.12.

\proclaim{\bf 1.14. Corollary {\rm \cite{Li2, Lemma 3.3 and Thms. 3.4, 3.9}}}
Either of the above conditions $(*)$ implies the following statements:
\smallskip

\noindent $a)$ The center $Z(G)$ of $G$ is contained in the
$\Cal H$-period subgroup $G_{\Cal H}$.
Moreover, the hypercenter
$Z_{\lim}(G) = \bigcup_{\alpha} Z_{\alpha}(G)$ of $G$
is contained in $G_{\Cal H}$, as well.
Hence, $X$ is Liouville whenever $G$ is hypernilpotent
or almost hypernilpotent\footnote
{Note that the condition $b)$ in \S 3.1 of \cite{Li2}
must sound as follows:
{\sl $[H,\, H_\alpha]\subseteq\bigcup_{\beta<\alpha} H_\beta$
for each $\alpha\in A$.}}.
\smallskip

\noindent $b)$ If the center $Z(G)$ of $G$ is nontrivial, then
the space $X$ is not Carath\'eodory hyperbolic. In particular,
a Galois covering $X$ over a quasiprojective variety cannot
be Carath\'eodory hyperbolic whenever its Galois group is an amenable
group with nontrivial center.
\smallskip

\noindent $c)$ If an element $s\in G$ is central in
a finite index subgroup $S\subseteq G$ $($or, more generally,
is contained in the centralizer of such a subgroup in $G)$,
then $s$ is an $\Cal H$-period: $s\in G_{\Cal H}$.
\endproclaim

\remark{\bf 1.15. Remarks} 1. The statement $Z(G)\subseteq G_{\Cal H}$
of Corollary 1.14$(a)$ follows
immediately from Theorem 1.13. Further, to prove the inclusion
$Z_{\lim}(G)\subseteq G_{\Cal H}$ one shows,
by transfinite induction, that the members
$Z_{\alpha}(G)$ of the transfinite upper central
series of $G$ are contained in $G_{\Cal H}$.
In turn, Corollary 1.14$(a)$ gives a proof of Theorem 1.6.
\smallskip

\noindent 2. We shall see in \S 3 that Corollary 1.14$(b)$
might be wrong, even for solvable coverings of compact
Riemann surfaces, if one omits the condition that the center
is nontrivial.
\smallskip

\noindent 3. Corollary 1.14$(c)$ shows that
if the conjugacy class
$s^G=\left\{g^{-1}sg\mid \,g\in G\right\}$ of an element $s\in G$
is finite, then any function $h\in \Cal H$ is constant on the $s^G$-orbit
$s^G x$ of any point $x\in X$. An element with the finite conjugacy class
is called an {\it $FC$-element}; in \S 2 we study some
generalizations and applications of this property.

On the other hand, in general, holomorphic functions need not be constant
on the orbits of infinite conjugacy classes, even for a free cocompact
holomorphic $G$-action of a (nonnilpotent) polycyclic group $G$
(see \S 4, especially Proposition 4.4$(b)$ and Remark 4.7).
\smallskip

\noindent 4. For some applications of Corollary 1.14$(b)$
see \cite{DetZa, DetOrZa}.
\endremark
\bigskip

\noindent{\bfit Some proofs}
\medskip

\noindent Following the scheme suggested in \cite{Li2}, we sketch here
the proofs of Proposition 1.10 and Theorem 1.13.
\vskip0.13pc

We denote by $\beta G$ the Stone-\v Cech compactification of
a discrete topological space $G$, or, which is the same,
the Gel'fand spectrum of the Banach
algebra $L^\infty (G)$ of all bounded complex valued functions on $G$.
Recall that the space $\beta G$ is compact and Hausdorff,
and $L^\infty (G) \cong C(\beta G)$.
For $f\in L^\infty (G)$ denote by $\hat f$ the unique continuous
extension of $f$ to $\beta G$, and by $M(f) \subseteq \beta G$
the peak point set of the function $\hat f$:
$$
M(f) =\left\{\xi\in\beta G\mid\,\modo{{\hat f}(\xi)}
=\normo{\hat f}_{C({\beta} G)}\right\}\,.
$$
The right action of a discrete group $G$
onto itself extends to the right $G$-action on $\beta G$.
\vskip0.13pc

Let $X$ be a Riemannian manifold resp. a complex space,
and let $G$ be a subgroup of the group $I(X)$ (see 1.5).
For any function $h\in {\Cal H}={\Cal H}(X)$
we set $\normo{h}_X = \sup_{x \in X} \modo{h(x)}$.
Let ${\Cal K}={\Cal K}(X)$ denote the convex cone of all nonnegative
bounded continuous subharmonic resp. plurisubharmonic functions on $X$.

\proclaim{1.16. Proposition {\rm \cite{Li2}}} Let $X$,\, $G$, and $\Cal H$
be as above. Assume that
\smallskip

\noindent $(i)$ the $G$-action on $X$ is ultra-Liouville,
i. e. the cone $\Cal K$ contains no nonconstant $G$-invariant function.
\smallskip

\noindent Let $h\in\Cal H$. Set $h_x(g) = h(gx)$ and $\varphi_h(x)
= \normo{\widehat{h_x}}_{C(\beta G)}$. Then
\vskip0.13pc

\noindent $a)$ $\varphi_h = \const$ \ and\qquad $b)$
$\normo{\widehat {h_x}}_{C({\beta} G)}
\equiv \varphi_h=\normo{h}_X\,$;
\vskip0.13pc

\noindent $c)$ the peak point set $M(\widehat {h_x}) \subseteq \beta G$
of the function $\widehat{h_x}$ does not depend on $x\in X$;
moreover, the subset $M(h):=M(\widehat {h_x})\subseteq\beta G$
is $G$-invariant;
\vskip0.13pc

\noindent $d)$ for any $G$-invariant regular probability
Borel measure $\mu$ on $\beta G$ the $L^2(\mu)$-class
$\left[\widehat{h_x}\right]$
of the function $\ \widehat{h_x}$ does not depend on $x\in X$.
\medskip

\noindent If, in addition, the group $G$ is amenable, then
\smallskip

\noindent $e)$ $\beta G$ carries a $G$-invariant probability
measure $\mu$ supported in $M(h)$;
\vskip0.13pc

\noindent $f)$ $h = 0$ whenever $\left[\widehat {h_x}\right] = 0$
in $L^2(\mu)$ for a measure $\mu$ as in $(e)$.
\endproclaim

\demo{Sketch of the proof} Note that the space
$\Cal H$ and the convex cone $\Cal K$
satisfy the following two conditions $(ii)$, $(iii)$:
\smallskip

\noindent $(ii)$ {\sl $\Cal H$ contains all the constant functions,
and for any $\normo{\cdot}_X$-bounded subset $\Cal F\subset \Cal H$
the function $k^2_{\Cal F}$, \
$k^2_{\Cal F}(x) = \sup_{f\in \Cal F}\modo{f(x)}^2$,
belongs to the cone $\Cal K$;}
\smallskip

\noindent $(iii)$ {\sl for any closed ball $B$
in the space $BC(X)$ of all complex valued bounded
continuous functions on $X$ the sets
${\Cal H}\cap B$ and ${\Cal K}\cap B$ are closed
in the compact open topology.}
\pagebreak 

Thus, by $(ii)$, \ $\varphi^2_h \in \Cal K$.
Since the function $\varphi^2_h$ is $G$-invariant,
$(a)$ follows from $(i)$. Clearly,
$$
\varphi_h\equiv \varphi_h(x)=\normo{\widehat{h_x}}_{C(\beta G)}
=\normo{h_x}_{L^\infty(G)}
= \sup_{g\in G} \modo{h(gx)}
\le \sup_{y\in X} \modo{h(y)}=\normo{h}_X\,.
$$
If the latter inequality were strict, then for some $x_\circ\in X$
we would have
$$
\aligned
\varphi_h<\modo{h(x_\circ)}
&\le \sup_{g\in G} \modo{h^g(x_\circ)}=\sup_{g\in G} \modo{h(gx_\circ)}\\
&=\sup_{g\in G} \modo{h_{x_\circ}(g)}
=\normo{\widehat{h_{x_\circ}}}_{C({\beta} G)}
=\varphi_h(x_\circ)=\varphi_h\,,
\endaligned
$$
which is impossible; this proves $(b)$.
\medskip

Given $x_\circ\in X$ and a point $\xi_\circ$
in the peak point set $M(\widehat{h_{x_\circ}})$,
consider the function $h^{\xi_\circ}(x)=\widehat{h_x}(\xi_\circ)$.
It follows from $(iii)$ that this function is in $\Cal H$,
and the function $\modo{h^{\xi_\circ}(x)}$
attains its maximal value $\normo{h}_X$ at the point
$x=x_\circ$. The maximum principle
\medskip

\noindent $(iv)$ {\sl $h = \const$ \ whenever \ $h\in \Cal H$ \
and \ $\modo{h(x_\circ)}=
\normo{h}_X$ for some point $x_\circ\in X$}
\medskip

\noindent implies that $h^{\xi_\circ}= \,$const. Hence,
$$
\modo{\widehat{h_x}(\xi_\circ)}=\modo{h^{\xi_\circ}(x)}
\equiv \modo{h^{\xi_\circ}(x_\circ)}
=\normo{h}_X=\modo{\widehat{h_{x_\circ}}(\xi_\circ)}\,.
$$
This shows that
$\xi_\circ\in M(\widehat{h_x})$ for any $x\in X$, which proves the first
assertion of $(c)$. The constant function $h^{\xi_\circ}$
is certainly $G$-invariant, and hence $h^{\xi_\circ g}=h^{\xi_\circ}$
for any $g\in G$. This yields
$\modo{\widehat{h_{x_\circ}}(\xi_\circ g)}
=\modo{h^{\xi_\circ g}(x_\circ)}
=\modo{h^{\xi_\circ}(x_\circ)}=\normo{h}_X$ and
$\xi_\circ g\in M(\widehat{h_{x_\circ}})= M(h)$,
which proves the second assertion of $(c)$.
\smallskip

Given a $G$-invariant regular probability Borel
measure $\mu$ on $\beta G$, define the function
$$
\Phi^2\colon \,X\to\Bbb R\,,\qquad
\Phi^2(x)=\normo{\left[\widehat{h_x}\right]}^2_{L^2(\mu)}
=\int_{\beta G} \modo{\widehat{h_x}(\xi)}^2 d\mu(\xi)\,.
$$
This function is $G$-invariant, and it follows from $(ii)$, $(iii)$
that $\Phi^2\in\Cal K$; \ by $(i)$, \ $\Phi^2= \,$const.
Fix a point $x_\circ\in X$, and consider
the mapping
$X\ni x\longmapsto F(x)=\left[\widehat{h_x}\right]\in L^2(\mu)$ and
the inner product $\psi(x)=\langle F(x)\,, F(x_\circ)\rangle$.
It follows from $(iii)$ that $\psi\in\Cal H$.
Clearly,
$$
\modo{\psi(x)}\le \normo{F(x)}_{L^2(\mu)} \normo{F(x_\circ)}_{L^2(\mu)}
= \Phi(x) \Phi(x_\circ)\equiv \Phi^2(x_\circ) \ \ \ {\text{and}} \ \ \
\modo{\psi(x_\circ)}=\Phi^2(x_\circ)\,;
$$
hence, by maximum principle $(iv)$, \,\,
$\langle F(x)\,, F(x_\circ)\rangle\equiv {\const}$.
Set $a=F(x_\circ)$ and $b=F(x)$; then we have
$\langle b,\, a\rangle=\normo{a}^2$ and $\normo{b}=\normo{a}$.
Since the norm in the Hilbert space $L^2(\mu)$ is strictly
convex, this implies $b=a$, that is, $F= \,$const. This proves $(d)$.
\smallskip

The statement $(e)$ follows from $(c)$ and the
Fixed Point Theorem for amenable groups \cite{Gre, Thm. 3.3.5}
applied to the natural $G$-action on the convex compact set of all
probability measures supported in the $G$-invariant set $M(h)$.
\smallskip

Finally, $(f)$ follows from
$(e)$. Indeed, the function $\widehat{h_x}$ is continuous on
$C(\beta G)$, and
hence $\left[\widehat {h_x}\right] = 0$
implies \ $\widehat{h_x}\mid {\supp}\,\mu = 0$.
Since ${\supp}\,\mu\subseteq M(h)=M(\widehat{h_x})$, it follows that
$\widehat{h_x}=0$ for any $x\in X$. Thus, $h=0$. 
\hfill $\square$
\enddemo
\pagebreak 

\remark{\bf 1.17. Remark} Let $X$ be a topological space endowed with a
$G$-action preserving a subspace ${\Cal H} \subset BC(X)$ and a convex cone
${\Cal K} \subset BC_{\Bbb R}(X)$. Suppose that the conditions $(i) - (iv)$
introduced above are fulfilled. {\sl All the assertions of
Proposition $1.16$ hold true in this more general setting}.
This yields analogs of Theorems 1.6 and 1.13
for certain equivariant second order elliptic operators
on smooth manifolds and for harmonic functions on discrete groups
\cite{Li2, \S\S 2.12--2.15}. Moreover, there is a version
of Proposition $1.16$ which applies to the case when, instead
of a $G$-action on $X$, one deals with $G$-actions in the space
$\Cal H$ and in the cone $\Cal K$. This leads to an analog of
Theorem 1.6 for suitable complex Lie group actions on
complex spaces (see \cite{Li2, \S2 and Thm. 2.17}).
\endremark
\medskip

\noindent {\bf 1.18.} {\it Proof of Proposition} 1.10.
$(a)$ is an immediate consequence of Proposition 1.16$(a, b)$.
To prove $(b)$, fix a function $h \in {\Cal H} (X)$ and define
the function $\widetilde h \in {\Cal H}(X \times X)$ by
$\widetilde h (x,\,y) = h(x) - h(y)$.
Clearly, $\widetilde h$
vanishes on the diagonal $\Delta \subset X \times X$, which 
is invariant under the diagonal $G$-action in $X\times X$;
hence, $\widetilde h$ vanishes on the orbit of any point
$(x,\, x)\in\Delta$. 
Since the diagonal action on $X \times X$ is assumed 
to be ultra-Liouville, the statement $(a)$ implies $\widetilde h = 0$.
Thus, $h = \const$ and $X$ is Liouville.
\hfill $\square$
\medskip

\noindent {\bf 1.19.} {\it Proof of Theorem $1.13$ for amenable $G$
and an element $s\in G$.} Actually, as in the proof of Proposition 1.16,
the only essential assumptions about the space $\Cal H$ and the convex
cone $\Cal K$ are those $(i)$-$(iv)$ above.
We deal with a function $f\in {\Cal H}$
and an element $s\in G$
such that $f^{[s,\,g]}=f$ for all $g\in G$, 
where $f^{[s,\,g]}(x) = f([s,\,g]x)=f(sgs^{-1}g^{-1}x)$. Thus,
$
f_{sx}(g)=f(gsx)=f(sgx)=f^s(gx)=(f^s)_x(g)
$
for all $g\in G$ and all $x\in X$, and hence
$$
\widehat{f_{sx}} = \widehat{(f^s)_x}\qquad{\text{for \ all}} \ x\in X\,.
\eqnum{*}
$$
Set $h=f^s-f\in{\Cal H}$. We must show that $h=0$.
Let $\mu$ be a measure as in Proposition 1.16$(e)$.
Since it is $G$-invariant, Proposition 1.16$(d)$ implies
that the $L^2(\mu)$-class $\left[\widehat{f_x}\right]$
does not depend on $x \in X$.
In particular, $\left[\widehat{f_{sx}}\right]=\left[\widehat{f_{x}}\right]$
and
$$
\left[\widehat{f_{sx}-f_x}\right]
=\left[\widehat{f_{sx}}-\widehat{f_{x}}\right]
=\left[\widehat{f_{sx}}\right]-\left[\widehat{f_{x}}\right]=0\,.
$$
Combined with $(*)$ this leads to
$\left[\widehat{ (f^s)_x-f_x} \right]=0$. Therefore,
$\left[\widehat{h_x}\right]
=\left[\widehat{\left( (f^s-f)_x \right)} \right]=0$.
Proposition 1.16$(f)$ implies $h=0$.
\hfill $\square$
\medskip

In the case when the element $s$ of $I(X)$ is not in $G$ and
there is no amenable subgroup in $I(X)$ containing
both $s$ and $G$, the above argument does not work.
To treat this case one should deal with actions of the amenable group
$G\times\Bbb Z$ in suitable function spaces
(see \cite{Li2} for details; see also Remark 1.17 above).
\smallskip

The proof of Theorem 1.13 for cocompact actions
is based on the compactness principle and a version of
the Harnack inequality. This approach goes back to Dynkin,
Malyutov \cite{DyMal} and Margulis \cite{Mar} who considered bounded
resp. positive harmonic functions on nilpotent groups
(see \cite{Li2} for details).


\head{\S 2. Upper FC-series and Liouville-type properties}\endhead

\noindent The concept of the upper FC-series \cite{Hai}
(see also Definition 2.1 below) allows us to generalize
Theorem 1.13 and Corollary 1.14.
Namely, let $X\to Y$ be a Galois covering with Galois group $G$
over an ultra-Liouville base $Y$.
We show that the period subgroup $G_{\Cal H}$
contains all the members of the upper FC-series of $G$
and, hence, their union, too  (see Corollary 2.5).
This leads also to a generalization of Theorem 1.6 and Corollary 1.8 on
the Liouville property of coverings (see Corollary 2.6).

\definition{2.1. Definitions} 1. {\it FC-groups and FC-series.} A group
$G$ is called {\it FC-group} \cite{Ba,\, Ku,\, To}
if the conjugacy class of each element $g\in G$ is finite.
Almost abelian groups or groups with finite
commutator subgroups are so \cite{Neu}.
Both of the latter classes contain the
proper subclass of groups $G$ with finite quotients $G/Z(G)$
by the center \cite{Neu;\, Er;\, To, Thm. 1.1}.
For any FC-group $G$ the quotient $G/Z(G)$ is a periodic group
(see \cite{Ba; To, Thm. 1.4}).

For any group $G$ the union $FC(G)$ of all finite conjugacy
classes is a normal subgroup of $G$. Clearly, $FC(G)$ is an FC-group;
it is called the {\it FC-center} of $G$ \cite{To}.
Set $FC_1(G) = FC(G)$ and for any $n\ge 1$ denote by $FC_{n+1}(G)$
the total preimage of $FC\left(G/FC_n(G)\right)$ under
the natural surjection $G\to G/FC_n(G)$. We obtain
the {\it upper FC-series} of $G$ \cite{Hai}:
$$
{\bold  1}  \vartriangleleft FC_1(G) \vartriangleleft FC_2(G)
\vartriangleleft \dots
\vartriangleleft FC_n(G) \vartriangleleft \dots \vartriangleleft G \,.
$$
Clearly, $FC_n(G)$ is a normal subgroup of $G$;
in fact, it is strictly characteristic\footnote
{A subgroup $H\subseteq G$ is called {\it strictly characteristic} if
$\phi(H)\subseteq H$ for any epimorphism
$\phi\colon \,G\to G$.}
\cite{Hai}. The upper FC-series extends transfinitely in the usual way
\cite{Du}, by defining $FC_\alpha(G)=\bigcup_{\beta<\alpha} FC_\beta(G)$
for each limit ordinal $\alpha$. We set
$$
FC_{\omega}(G) = \bigcup_{n \in \Bbb N} FC_n(G)
\qquad\text{and}\qquad
FC_{\lim}(G)=\bigcup_{\alpha} FC_\alpha (G)\,,
$$
where $\alpha$ runs over all the ordinals. The normal
subgroup $FC_{\lim}(G)\vartriangleleft G$ is called the
{\it FC-hypercenter} of $G$; we say that the elements of
$FC_{\lim}(G)\vartriangleleft G$ are {\it FC-hypercentral} in $G$.
\vskip0.13pc

\noindent 2. {\it FC-nilpotent and FC-hypernilpotent groups} \cite{Hai, Du, Rob}.
If $G = FC_n(G)$ for some $n\in\Bbb N$ and $G \ne FC_{n-1}(G)$,
then $G$ is called {\it FC-nilpotent of class $n$},
or simply {\it FC-nilpotent}.
We say that the group $G$ is
{\it FC-$\omega$-nilpotent} resp. {\it FC-hypernilpotent}
if $G = FC_{\omega}(G)$ resp. $G = FC_{\lim}(G)$.
Clearly, an FC-$\omega$-nilpotent group is {\it locally FC-nilpotent},
meaning that any finitely generated subgroup of $G$ is FC-nilpotent.
\enddefinition

\remark{\bf 2.2. Remarks} 1. A nilpotent (resp. $\omega$-nilpotent,
hypernilpotent, locally nilpotent) group is FC-nilpotent
(resp. FC-$\omega$-nilpotent, FC-hypernilpotent, locally FC-nilpotent).
The properties of FC-nilpotence,
FC-$\omega$-nilpotence,
and FC-hypernilpotence are inherited by the subgroups, the quotient groups,
and the finite extensions.
So, a finite
extension of a nilpotent group
is FC-nilpotent. Vice versa, a finitely generated FC-nilpotent group of
class $n$ is a finite extension of a nilpotent
group of class at most $n$ (see \cite{DuMcL, Thm. 2}).
{\sl A finitely generated FC-hypernilpotent
group is almost nilpotent and almost
polycyclic} \cite{McL; Rob, Vol. 1, p. 133}. Thus, an
FC-hypernilpotent group is locally FC-nilpotent.

\noindent 2. Every locally FC-nilpotent, and so, every FC-hypernilpotent,
group $G$ is amenable. Indeed, $G$ is the union of the direct system of
its finitely generated FC-nilpotent subgroups.
Therefore, by Theorem 1.2.7 in \cite{Gre},
the statement follows once we know that any
finitely generated FC-nilpotent group is amenable. The latter holds
since a finitely generated FC-group is almost nilpotent
\cite{DuMcL, Thm. 2} (see Remark 2.2.1 above).
\smallskip

\noindent 3. The following examples show that, in general,
the FC-(hyper)center of a finitely generated group need not be
finitely generated, and a countable FC-group which is not
finitely generated may be neither almost solvable nor almost
hypernilpotent.
\endremark
\medskip

\example {2.3. Examples} 1 (see \cite{PHa} and also
\cite{Rob, Vol. 1, Thm. 5.36}). Let
$$
G=\langle a,b\mid \,\left[b_i,\,b_j,\,b_k\right]=1, \
\left[b_i,\,b_j\right]=\left[b_{i+k},\,b_{j+k}\right], \
i,\,j,\,k\in\Bbb Z\rangle\,,
$$
where $b_i=a^{-i}ba^i$. Then the center $Z(G)$ coincides
with the subgroup $H\subset G$ generated by the elements
$d_r=[b_0,\, b_r]$, \ $r\in\Bbb N$. It is a free abelian
group of infinite rank. Furthermore, $Z(G/H)=\bold 1$,
and hence $Z_{\lim}(G)=Z(G)=H$. It is easily seen that
actually $FC_{\lim}(G)=FC(G)=Z(G)=H$.
\medskip

\noindent 2 (see \cite{Ku, \S 38}). Let
$$
P=\underset{n=2k+1,\, k\ge 2}\to\prod A_n\qquad\text{and}\qquad
H=\underset{n=2k+1,\, k\ge 2}\to\bigoplus A_n
$$
be, respectively, the direct product
and the direct sum\footnote
{i. e. the restricted direct product.}
of the alternating groups $A_n\subset S_n$
of all odd degrees $n\ge 5$, where $S_n$ stays for the symmetric group
of degree $n$.
It is known that for any odd $n$ the group $A_n$
is generated by the cyclic permutations
$a_n=(1,...,n)$ and $b_n=(1,2,3)$.
Let $G\subset P$ be the subgroup generated
by the elements $a=(a_5,...,a_{2k+1},...)\in P$ and
$b=(b_5,...,b_{2k+1},...)\in P$. Then $G\supset H$.
It is easily seen that the center $Z(G)$ is trivial,
whereas the FC-center $FC(G)$ coincides with the subgroup $H$.
Moreover, $FC(G/H)=\bold 1$, and hence
$FC_{\lim}(G)=FC(G)=H$. The FC-group $H=FC(G)$ is neither
almost solvable nor almost hypernilpotent. Indeed,
for any normal subgroup $N$ of finite index in $H$
the intersection $N_n=A_n\cap N \ne \bold 1$ if $n$ is sufficiently large.
Clearly, $A_n=N_n\subset N$ for such $n$ (since $A_n$ is simple), and thus
$N$ is neither solvable nor hypernilpotent.
\endexample
\medskip

The following lemma is an easy consequence of Theorem 1.13.
We use the same notation
as in Theorem 1.13; see also Definitions 1.5, 1.12.

\proclaim{2.4. Lemma} Suppose that one of the
two conditions of Theorem $1.13$ is fulfilled, i. e. either
\smallskip

\noindent $\,*$ $G$ is amenable and its action on $X$ is ultra-Liouville, or
\smallskip

\noindent $\,*$ the action of $G$ on $X$ is cocompact.
\smallskip

\noindent Let $N \vartriangleleft G$ be a normal subgroup of $G$ contained
in the $\Cal H$-period subgroup $G_{\Cal H}$.
Let $\bar s$ denote the image of an element $s\in G$ in the
quotient group $\overline {G}=G/N$. Suppose that the conjugacy class
$\bar {s}^{\overline G}
=\left\{\bar {g}^{-1}\bar {s}\bar {g}
\mid \,\bar {g}\in\overline{G}\right\}$ of $\bar s$
in $\overline G$ is finite. Then $s\in G_{\Cal H}$.
\endproclaim
\pagebreak 

\demo{Proof} By our assumption, the centralizer ${\overline C}$
of the element $\bar s$ is of finite index in $\overline G$.
The total preimage $C$ of $\overline C$ in $G$
is a subgroup of finite index.
Hence (see \cite{Li2,\, Lemma 3.3}) $C$ satisfies the same
condition ($*$) as $G$. Furthermore, $C$
contains both $s$ and $N$. Since $\bar s$ is central in
$\overline C$, we have $[s,\,C]\subseteq N\subseteq G_{\Cal H}$. Now
Theorem 1.13 shows that $s\in C_{\Cal H}\subseteq G_{\Cal H}$.
\hfill $\square$
\enddemo

\proclaim{2.5. Corollary\footnote{cf. Corollary 1.14$(a)$.}}
Suppose that one of the conditions
$(\,*)$ of Lemma $2.4$ is fulfilled. Then
$FC_{\lim}(G) \subseteq G_{\Cal H}$, i. e. any FC-hypercentral
element of $G$ is an ${\Cal H}$-period.
\endproclaim

\demo{Proof} Starting with the unit subgroup
${\bold 1}\subseteq G_{\Cal H}$, we proceed by transfinite induction.
Suppose that $FC_{\alpha}(G)\subseteq G_{\Cal H}$.
Set $N = FC_{\alpha}(G)\vartriangleleft G$. By Lemma 2.5,
for any element $s\in FC_{\alpha+1}(G)$ we have $s\in G_{\Cal H}$,
and thus $FC_{\alpha+1}(G)\subseteq G_{\Cal H}$.
Furthermore, if $\alpha$ is a limit ordinal
and $FC_{\beta}(G)\subseteq G_{\Cal H}$ for all $\beta < \alpha$,
then $FC_{\alpha}(G)
= \bigcup_{\beta < \alpha} FC_{\beta}(G)\subseteq G_{\Cal H}$.
By induction, it follows that $FC_{\lim}(G)
= \bigcup_{\alpha} FC_{\alpha}(G)\subseteq G_{\Cal H}$.
\hfill $\square$
\enddemo

\proclaim{2.6. Corollary\footnote{cf. Theorem 1.4 and Corollary 1.8.}}
Let $X\to Y$ be a Galois covering with Galois group $G$ over
a compact Riemannian resp. K\"ahler manifold $Y$. If $G$ is an extension
of an FC-hypernilpotent group by a Varopoulos group,
then $X$ is Liouville.
\endproclaim

\remark{\bf 2.7. Remark} It follows from Corollary 2.5 that
either of the conditions ($*$) of Lemma 2.4 implies
the following property of the period subgroup
$G_{\Cal H}\vartriangleleft G$: \
the FC-center $FC(G/G_{\Cal H})$ of the quotient
$G/G_{\Cal H}$ is trivial, that is, each
nontrivial (i. e. $\ne \{\bold 1\}$) conjugacy class in $G/G_{\Cal H}$
is infinite.
Clearly, the subgroup $FC_{\lim}(G) \vartriangleleft G_{\Cal H}$
has the same property. It would be interesting to find
an example (if it does exist) in which  $FC_{\lim}(G) \ne G_{\Cal H}$.
\endremark


\head{\S 3. Solvable Carath\'eodory hyperbolic coverings
of a compact Riemann surface}\endhead

\noindent It was shown in \cite{LySu} that any
compact Riemann surface $Z$ of genus $g\ge 2$ admits
a non-Liouville Galois covering $X\to Z$
with a metabelian (i. e. two-step solvable)  Galois group.
Modifying the construction of
Lyons and Sullivan, we prove the following theorem.

\proclaim{3.1. Theorem} Each compact Riemann surface $Z$ of genus
$g\ge 2$ admits a Carath\'eodory hyperbolic metabelian covering
$X\to Z$.
\endproclaim

\demo{Proof} Let $\Gamma = \pi_1(Z)$ and let $p_1\colon\,Y\to Z$ 
be the maximal abelian covering over $Z$ 
(i. e., the covering corresponding to 
the commutator subgroup $\Gamma^\prime=[\Gamma,\,\Gamma]$);
its Galois group
$G= \Gamma/\Gamma^\prime\cong H_1(Z,\,\Bbb Z)\cong \Bbb Z^{2g}$.
Since $\rk\, G = 2g\ge 4$, a theorem of A. Mori \cite{Mo} implies that
for an any point $y\in Y$ there exists
a unique positive Green function, say $g_y$, with pole at $y$
(see also \cite{Ts, Theorem X.46}).

Let $D\subset Y\setminus \{y\}$ be a simply connected domain. Then
there is a conjugate harmonic function $g_y^*$ of $g_y$ in $D$, which
is defined uniquely up to an additive real constant.
Therefore, the differential $\omega_y = df_y$, where $f_y = g_y + i g_y^*$,
is a well-defined holomorphic $1$-form on $Y\setminus \{y\}$.
Its real part $\Re\,\omega_y = dg_y$ is an exact $1$-form
on $Y\setminus \{y\}$. Hence,
the real part of each period $\int_{\gamma}\omega_y$ of $\omega_y$, where
$\gamma\in H_1 (Y\setminus\{y\};\,\Bbb Z)$, is zero.
Thus, $\omega$ defines a homomorphism
$H_1 (Y \setminus \{y\};\,\Bbb Z) \to i\Bbb R$.
\smallskip

Fix a point $z_0 \in Y \setminus \{y\}$.
For any particular choice of $g_y^*$ consider the function
$$
\varphi_y (z) = \exp \, \left(-2\pi f_y(z)\right)
= \exp \, \left(-2\pi\left(f_y\left(z_0\right)
+ \int_{z_0}^z \omega_y\right)\right)\,.
$$
This is a multi-valued holomorphic function on $Y$ with values in
the unit disk $\Bbb D$. For a given $y \in Y$ any two such functions coincide
up to a constant factor $\lambda \in \Bbb T$, where
$\Bbb T=\left\{\lambda\in\Bbb C\mid \,\modo{\lambda} = 1\right\}$.
For each $y \in Y$ choose, once forever, one of the functions $\varphi_y$.
\smallskip

Any two values of $\varphi_y$ differ by a factor of the form
$\exp \, \left(-2\pi\int_{\gamma}\omega_y\right) \in \Bbb T$,
where $\gamma \in H_1 (Y;\,\Bbb Z)$. More precisely, we have a well-defined
character
$$
\alpha_y\colon \,H_1 \left(Y \setminus \{y\}; \, \Bbb Z\right) \ni
\gamma \longmapsto  \exp \, \left(-2\pi\int_{\gamma}\omega_y\right) \in \Bbb T\,.
$$
Actually, it yields a character
$$
\alpha_y\colon \,H_1 (Y ;\,\Bbb Z) \to \Bbb T\,.
$$
Indeed, consider the exact sequence
$$
0 \to \Bbb Z \to H_1 (Y \setminus \{y\};\,\Bbb Z)
\to H_1 (Y ;\,\Bbb Z) \to 0\,,
$$
where the subgroup $\Bbb Z \subset H_1 (Y \setminus \{y\} ;\,\Bbb Z)$
is generated by a small circle $\sigma_{\epsilon}$
in $Y$ centered at $y$.
In a small disk $\delta_y$ around $y$ we have
$g_y(z) = -{1\over 2\pi}\log\modo{z - y} + h_y(z)$, where
$h_y$ is a single-valued harmonic function in $\delta_y$; hence,
$f_y(z) = -{1\over 2\pi}\log(z - y) + \widetilde{f_y}(z)$, where
$\widetilde{f_y}$ is a single-valued holomorphic function in $\delta_y$.
It follows that
$$
2\pi \int_{\sigma_\epsilon} \omega_y = 2\pi \int_{\sigma_\epsilon}
\left (- {1\over 2\pi}d\log(z - y)
+ d\widetilde{f_y} \right)\in 2 \pi i \Bbb Z\,.
$$
Thereby, $\exp \, \left(-2\pi \int_{\sigma_y} \omega_y\right)=1$,
the restriction of the homomorphism $\alpha_y$
to the kernel subgroup $\Bbb Z$ in the above exact sequence is trivial,
and $\alpha_y$ can be pushed down to the quotient group.

The set of values of the function $\varphi_y$ at a point
$z \in Y \setminus \{y\}$ coinsides with a coset of the
subgroup Image$\,(\alpha_y)$ in the multiplicative group $\Bbb C^*$,
whereas all its values at the point $y$ are zero.

Let $\rho\colon \,\pi_1(Y) \to H_1(Y; \,\Bbb Z) \cong
\pi_1(Y)/\pi_1^\prime(Y)$ be  the canonical surjection.
Set $\widetilde {\alpha_y} = \alpha_y \circ \rho$.
The covering $X_y \to Y$ over $Y$ corresponding to the subgroup
$\Ker\, \widetilde {\alpha_y} \vartriangleleft \pi_1(Y)$ is the minimal one
such that the function $\varphi_y$ becomes single-valued when lifted
to $X_y$. Set
$$
K = \bigcap_{y \in Y} {\Ker}\,\alpha_y \subset H_1(Y; \,\Bbb Z)\,,
$$
and ${\widetilde K} = \rho^{-1}(K) \subset \pi_1(Y)$.
Let $p_2\colon \,X \to Y$ be the abelian
covering over $Y$ associated with the subgroup
${\widetilde K} \vartriangleleft \pi_1(Y)$. First we show that
$p=p_1\circ p_2\colon\,X \overset{{p_2}}\to{\longrightarrow} Y
\overset{{p_1}}\to{\longrightarrow} Z$ is a Galois
(and hence a metabelian) covering.

Let $\Gamma^{\prime\prime}=[\Gamma^\prime,\,\Gamma^\prime]$.
The group $G= \Gamma/\Gamma^\prime\cong H_1(Z,\,\Bbb Z)$
acts isometrically on $Y$; thus,
$g_y \circ \tilde\gamma = g_{\tilde\gamma(y)}$ and
$\tilde\gamma^*(\omega_y) = \omega_{\tilde\gamma(y)}$
for any $\tilde\gamma \in \Gamma/\Gamma^\prime$.
Hence the subgroup
$K\subset H_1(Y;\,\Bbb Z)\cong\Gamma^\prime/\Gamma^{\prime\prime}$
is invariant with respect to the induced action of $\Gamma/\Gamma^\prime$
in homology $H_1(Y;\,{\Bbb Z})$; denote this action
as $\mu$. It is easily seen that $\mu$ coincides with the ``adjoint"
representation
$\Gamma/\Gamma^\prime\ni\tilde\gamma
\mapsto T_{\tilde\gamma}\in\Aut\,(\Gamma'/\Gamma'')$, \
$T_{\tilde\gamma}(v) = {\tilde\gamma}^{-1}v{\tilde\gamma}$, \
$v\in \Gamma^\prime/\Gamma^{\prime\prime}$.
It follows that
$p_{*}(\pi_1(X)) = (p_1)_{*}({\widetilde K}) \subset \pi_1(Z)$ is a normal
subgroup of $\pi_1(Z)$, and hence $p\colon\,X\to Z$
is a metabelian Galois covering.
\medskip

Clearly, $p_2\colon\,X\to Y$ is the minimal covering over $Y$ such
that all the functions
$\{\varphi_y\}_{y \in Y}$ become single-valued when lifted to $X$. Let
$E = \left \{\widetilde {\varphi_y}\right \}_{y \in Y} \subset H^{\infty}(X)$
be the collection of all the lifted functions.
We will show that $E$ separates the points of $X$.

Denote by $F_y = p^{-1}(y) \subset X$ the fiber of $p$
over $y \in Y$. For any two
distinct points $y, y' \in Y$ the function $\widetilde {\varphi_y}$ vanishes
identically on $F_y$ and does not vanish at the points of  $F_{y'}$.
Therefore, $E$ separates the fibers $\{F_y\}$.

Thus, it is sufficient to show that $E$ separates the
points of each fiber $F_y$. It is easily seen that for $y' \ne y$
the function $\widetilde {\varphi_y}$ separates the points of $F_{y'}$ if and
only if
${\Ker}\,\alpha_y = K$. If the latter equality holds
for a certain pair of distinct
points $y_1,\,y_2 \in Y$, then the points of each fiber
$F_y$, \ $y \in Y$, are separated by at least one of the functions
$\widetilde {\varphi_{y_1}}$, \ $\widetilde {\varphi_{y_2}}$.
Hence, the theorem follows from the next claim.
\medskip

\noindent {\bf Claim 1.} {\sl There exists a countable union
\ ${\Cal C} = \bigcup_{n \in \Bbb N} C_n \subset Y$
of real analytic curves $C_n$ in $Y$ such that
${\Ker}\,\alpha_y = K$ for each point $y \in Y \setminus {\Cal C}$.}
\medskip

The proof is based on the following statement\footnote{
It should be well known; for the sake of completeness we
give a simple proof.}:
\medskip

\noindent {\bf Claim 2.} {\sl The function of two complex variables
$g(y,\,y') = g_y(y')$ is harmonic on the complement
$(Y\times Y) \setminus \Delta$,
where $\Delta \subset Y\times Y$ is the diagonal.}
\medskip

\noindent {\it Proof of Claim $2$.} By the symmetry property of Green function
\cite{Ts, Theorem I.16}, we have $g_y(y') = g_{y'}(y)$ for any
$y \ne y'$, \ $y,\,y' \in Y$.  Hence, $g(y,\,y')$
is a harmonic function in each argument on
$(Y\times Y) \setminus \Delta$.
It is sufficient to show that it is harmonic as a function of two
complex variables in each bidisk
$\delta \times \delta'\subset\subset (Y\times Y)\setminus \Delta$,
where $\delta$, \ $\delta'$ are two small disks in $Y$.
Being harmonic in each variable,
the function $g(y,\,y')$ in the bidisk $\delta \times \delta'$
satisfies the Laplace equation
$\Delta_{y,\, y'} g(y,\,y')
=\Delta_y g(y,\,y') + \Delta_{y'} g(y,\,y') =0$,
where $\Delta_y$, $\Delta_{y'}$ are the usual Laplacians.
Therefore, $g(y,\,y')$ is harmonic in $\delta \times \delta'$
as soon as it is continuous there.

Since the function $g(y,\,0)$ is continuous in the closed disk
$\overline{\delta}$, the family $g_y=g_y(y')$ of positive harmonic
functions in $\delta'$ is equicontinuous in every smaller
closed disk (this easily follows by the Harnack
inequality). This implies that $g=g(y,\,y')$ is a continuous function in
$\delta \times \delta'$, which completes the proof. \hfill $\square$
\medskip

\noindent {\it Proof of Claim $1$.} It is sufficient to check our
statement locally. Fix a small disk $\delta \subset Y$. We will
show that ${\Ker}\,\alpha_y = K$ for all $y \in \delta$
outside of a countable union
${\Cal C}_{\delta} \subset \delta$ of closed real analytic curves in $\delta$.

It follows from Claim 2 that in each local chart $\Omega$ in $Y$
the coefficients of the holomorphic $1$-form $\omega_y$ are real
analytic functions of $y \in Y \setminus \Omega$.

Let a sequence $\left\{\gamma_n\right\}_{n \in \Bbb N}$ of $1$-cycles in $Y$
be a free basis of the homology group
$$
H_1(Y; \,\Bbb Z) \cong \Bbb Z^{\infty}
= \bigoplus\limits_{1}^{\infty} \Bbb Z\,.
$$
We may assume that they do not meet the closed disk
$\overline \delta$. 
The periods $c_n(y) = \int_{\gamma_n} \omega_y$, \ $n \in \Bbb N$,
are (imaginary-valued) real analytic functions of $y \in \delta$.. For
$\gamma = \sum_{j=1}^n  a_j\gamma_j \in H_1(Y; \,\Bbb Z)$ we have
$$
\langle \gamma,\,\omega_y \rangle = \int_{\gamma} \omega_y
= \sum_{j=1}^n  a_j\int_{\gamma_j} \omega_y = \sum_{j=1}^n  a_jc_j(y) =
\langle {\overline a},\,c(y) \rangle\,,
$$
where ${\overline a} = (a_1,\,\dots,a_n,\,0,\dots)$ and
$c(y) = (c_j(y))_{j=1}^{\infty}$. By the definition
of the character $\alpha_y\colon \,H_1(Y; \,\Bbb Z)\to \Bbb T$, we have
$$
\Ker\,\alpha_y
= \left\{\gamma\in H_1(Y; \,\Bbb Z)\mid \,\langle \gamma,\, i \omega_y \rangle
= \langle {\overline a},\, i c(y) \rangle \in \Bbb Z\right\}\,.
$$
Set
$L = \left\{{\overline a}\in\Bbb Z^{\infty}\mid
\,\langle {\overline a},\, i c(y)\rangle
\in \Bbb Z  \,\,\,{\text{for \,\,all}}\,\,\,y \in \delta\right\}$.
For each ${\overline a} \in \Bbb Z^{\infty} \setminus L$
and for each $k \in \Bbb Z$
consider the real analytic curve
$C_{\overline a,\,k}
= \left\{y\in\delta\mid \,\langle {\overline a},\, i c(y) \rangle
= k\right\}$.
Put
$$
{\Cal C}_{\delta} = \bigcup_{{\overline a}\in \Bbb Z^{\infty}\setminus L;\,\,\,
k \in \Bbb Z} C_{\overline a,\,k}\,.
$$
It is easily seen that for $y \in \delta \setminus {\Cal C}_{\delta}$ the
subgroup
${\Ker}\,\alpha_y \subset H_1(Y; \,\Bbb Z)$ does not depend on $y$
and coincides with $L$. Furthermore, for any
$y \in {\Cal C}_{\delta}$ we have ${\Ker}\,\alpha_y \supset L$,
and hence $L = K$.
This proves Claim 1 and completes the proof of Theorem 3.1.
\hfill $\square$
\enddemo

\head{\S 4. $H^\infty$-hulls in a solvable covering of Inoue surface}\endhead

\noindent In this section we study in more details
the universal covering $\pi\colon \,X\to {\Cal I}$
over one of the Inoue surfaces ${\Cal I}$ \cite{In}.
We start with a description of the Inoue surface.
\medskip

Let $A\in\SL\,(3;\,\Bbb Z)$ be a matrix with one real
eigenvalue $\alpha>1$ and two complex conjugate eigenvalues
${\beta},\,\overline{\beta} \in \Bbb C \setminus \Bbb R$ (certainly
$\modo{\beta}<1$).
Let ${\bold  a}=(a_1,\, a_2,\, a_3)$ resp.
${\bold  b}=(b_1,\, b_2,\, b_3)$ be
a real resp. a complex eigenvector of $A$ corresponding to the
eigenvalue $\alpha$ resp. $\beta$.

Set $\Bbb H=\{z\in\Bbb C\mid \,\Im\,z>0\}$ (the upper halfplane) and
$X=\Bbb H\times \Bbb C$. Consider the subgroup $G\subset \Aut\, X$
generated by the following four automorphisms $g_j$:
$$
g_0(z,w)=(\alpha z,{\beta} w),\qquad
g_j(z,w)=(z+a_j, w+b_j), \ \ 1\le j\le 3, \ \
(z,w)\in X=\Bbb H\times \Bbb C\,.
$$
The action of the group $G$ on $X$ is free,
properly discontinuous, and cocompact. The smooth compact complex
surface ${\Cal I}=X/G$ is one of the {\it Inoue surfaces} \cite{In}.

The subgroup $G_0\subset G$ generated by $g_1,\, g_2,\, g_3$ is isomorphic
to $\Bbb Z^3$; this subgroup is normal in $G$,
and the quotient group $G/G_0$ is isomorphic to $\Bbb Z$.
Thus, we have the exact sequence
$$
0\longrightarrow \Bbb Z^3 \longrightarrow G
{\overset{\tau}\to{\longrightarrow}} \,\Bbb Z \longrightarrow 0\, ,
\eqnum{1}
$$
and the corresponding tower of abelian coverings
$X \ {\overset{\Bbb Z^3}\to{\longrightarrow}} \ Y
{\overset{\Bbb Z}\to{\longrightarrow}} {\Cal I}$,
where $Y=X/G_0$.
In particular, $G$ is a metabelian (i. e. two-step solvable)
polycyclic group, and $X\to \Cal I$ is a polycyclic covering
with the Galois group $G$.
\medskip

To establish certain analytic properties
of the covering $X\to\Cal I$,
we need the following simple observations.

First, note that the sequence (1) splits;
a splitting $\rho\colon \,\Bbb Z\to G$
\ ($\tau\circ\rho=\id_{\Bbb Z}$) may be defined by
$\Bbb Z\ni m\mapsto g_0^m\in G$. Therefore, $G$ is a semidirect product
$\Bbb Z^3 \leftthreetimes\Bbb Z$, and any element $g\in G$
admits a unique representation of the form
$$
g=g_0^m\widetilde g=g_0^m g_1^{r_1}g_2^{r_2}g_3^{r_3},
\,\,\,\,{\text{where}} \,\,\,\,
m=\tau(g)\in\Bbb Z, \, \, r_1,r_2,r_3\in\Bbb Z, \,\,\,
\widetilde g=g_1^{r_1}g_2^{r_2}g_3^{r_3}\in G_0.
$$
Using this normal form, for any $d\in \Bbb Z$ we can write
$$
\aligned
g^{-1}g_0^d g &=(g_0^m\widetilde g)^{-1}
	    \cdot g_0^d\cdot(g_0^m\widetilde g)\\
&\\
&=\widetilde g^{-1}g_0^{-m}\cdot g_0^d\cdot
		     g_0^m\widetilde g=\widetilde g^{-1}g_0^d\widetilde g
		     =g_3^{-r_3}g_2^{-r_2}g_1^{-r_1}\cdot g_0^d\cdot
		      g_1^{r_1}g_2^{r_2}g_3^{r_3}.
\endaligned
$$

\proclaim{4.1. Lemma} The conjugacy class $s^G$ of the element
$s=g_0^d$ consists of all the transformations of the form
$$
(z,w)\mapsto\left (\alpha^d z + (\alpha^d-1)(r_1 a_1 + r_2 a_2 + r_3 a_3),
\, {\beta}^d w + ({\beta}^d-1)(r_1 b_1 + r_2 b_2 + r_3 b_3)\right )\,,
$$
where \ $r_1$, $r_2$, and $r_3$ run over $\Bbb Z$.
\endproclaim

\demo{Proof} Since the elements $g_j$, \ $j=1,2,3$, commute,
the lemma follows from (2) and the formula
$$
\ \ \qquad\qquad\qquad g_j^{-r}g_0^d g_j^r(z,w)
=(\alpha^d z + (\alpha^d-1)ra_j, {\beta}^d w + ({\beta}^d-1)rb_j).
\qquad\qquad\qquad \ \ \square
$$
\enddemo

\proclaim{4.2. Lemma} $a)$ The real eigenvalue
$\alpha$ of the matrix $A$ is a nonquadratic irrationality.
\smallskip

\noindent $b)$ The coordinates $a_1, a_2, a_3$ of the corresponding eigenvector
$\bold  a$ are linearly independent over $\Bbb Q$.
\smallskip

\noindent $c)$ For any subgroup $L\subseteq\Bbb Z^3$ of rank rk$\, L\ge 2$
and for any finite subset $S\subset L$ we have
$$
\inf\left\{\modo{r_1 a_1 + r_2 a_2 + r_3 a_3}\mid \,
{\bold  r}=(r_1,\, r_2,\, r_3)\in L\setminus S\right\} = 0.
$$
\endproclaim

\demo{Proof} $a)$ The characteristic polynomial
$P(z)=t^3+pt^2+qt+1$, \ $p,\,q \in \Bbb Z$,
of the unimodular matrix $A$ has no rational root
except, possibly, of $\pm 1$. Since $\alpha,\,{\beta}\ne\pm 1$, the
polynomial $P$ is irreducible over $\Bbb Q$, which proves $(a)$.
\smallskip

\noindent $b)$ Assume, on the contrary, that $a_1,\,a_2,\,a_3$
are linearly dependent
over $\Bbb Q$. Let $A=(a_{ij})_{i,j=1}^3$, where $a_{ij}\in \Bbb Z$.
Then we have the following system
$$
\aligned
r_1 a_1 &+\qquad \ \ \ r_2 a_2 +\qquad \ \ \ \ r_3 a_3=0\\
(a_{11}-\alpha)a_1&+\qquad \ \ a_{12}a_2+\qquad \ \ a_{13}a_3=0\\
a_{21}a_1&+(a_{22}-\alpha)a_2+\qquad \ \ a_{23}a_3=0\\
a_{31}a_1&+\qquad \ \ a_{32}a_2+(a_{33}-\alpha)a_3=0\\
\endaligned
$$
with some $(r_1,\,r_2,\,r_3)\in\Bbb Z^3\setminus\{\bold  0\}$.
Since $(a_1,\,a_2,\,a_3)\ne\bold 0$, we obtain the following three equations
(each of degree at most $2$) for $\alpha$:
$$
\det\pmatrix
r_1 & r_2 & r_3\\
a_{21} & a_{22}-\alpha & a_{23}\\
a_{31} & a_{32} & a_{33}-\alpha
\endpmatrix
=0\,,\qquad
\det\pmatrix
a_{11}-\alpha & a_{12} & a_{13}\\
r_1 & r_2 & r_3\\
a_{31} & a_{32} & a_{33}-\alpha
\endpmatrix
=0\,,
$$
$$
\det\pmatrix
a_{11}-\alpha & a_{12} & a_{13}\\
a_{21} & a_{22}-\alpha & a_{23}\\
r_1 & r_2 & r_3
\endpmatrix
=0\,.
$$
At least one of these equations must certainly be of degree $2$
(for $(r_1,\,r_2,\,r_3)\ne\bold  0$), which contradicts $(a)$.
\smallskip

\noindent $c)$ By $(b)$, the homomorphism
$
\chi\colon L\ni{\bold  r}=(r_1,\,r_2,\,r_3)\mapsto
r_1 a_1 + r_2 a_2 + r_3 a_3\in\Bbb R
$
is injective. Hence, $M=\chi(L)\subset\Bbb R$ is a free abelian
subgroup of rank $\rk\, M = \rk\, L\ge 2$. The closure
$\overline M$ of $M$ in $\Bbb R$ coincides with $\Bbb R$ (for otherwise,
$\overline M\cong\Bbb Z$ and hence $M\cong\Bbb Z$, which
contradicts the property $\rk\, M\ge 2$). This implies $(c)$.
\hfill $\square$
\enddemo

\definition{4.3. Definition} The $H^\infty(X)$-{\it hull} \
$\widehat Y$
of a set $Y\subseteq X$ in a complex space $X$ is defined as follows:
$$
\widehat Y = H^\infty{-\hull}_{X}(Y) =
\left\{x\in X\mid \,\modo{f(x)}\le \sup_{y\in Y} \modo{f(y)}
\ \ {\text{for \ all}} \ \ f\in H^\infty(X)\right\}\,.
$$
\enddefinition

\proclaim{4.4. Proposition {\rm (cf. Corollary 1.14($c$))}} Let,
as in Lemma $4.1$, \ $s=g_0^d$. Suppose that $\alpha^d > 2$. Then:
\smallskip

\noindent $a)$ The $s^G$-orbit $s^G (x_\circ)
=\{g^{-1}sgx_\circ\mid \,g\in G\}$
of the point $x_\circ=(i,0)\in X=\Bbb H\times\Bbb C$
consists of all the points \ $x = (z,w) \in X=\Bbb H\times\Bbb C$ \
of the form
$$
(z,w) = \left ((\alpha^d-1)
(r_1 a_1 + r_2 a_2 + r_3 a_3) + i\alpha^d, \ \
({\beta}^d-1)(r_1 b_1 + r_2 b_2 + r_3 b_3)\right ),
$$
where $r_1,\, r_2,\, r_3\in\Bbb Z$.
\smallskip

\noindent $b)$ The bounded holomorphic function
$F(z,w) =
2(z+i)^{-1}$ on $X=\Bbb H\times\Bbb C$ satisfies the inequality
$$
\modo{F(x_\circ)}=1 > {2\over 3} > \sup_{x\in s^G (x_\circ)}\modo{F(x)}.
$$
In particular, the $H^\infty(X)$-hull $\widehat {{s^G(x_\circ)}}$
of the $s^G$-orbit $s^G(x_\circ)$ does not contain the
point $x_\circ$ itself,
and for any mean $\frak m$ on $L^\infty(s^G(x_\circ))$ we have
$F(x_\circ)\ne\frak m (F\mid{s^G(x_\circ)})$.
\endproclaim

\demo{Proof} $(a)$ follows immediately from Lemma 4.1.
In view of the assumption $\alpha^d > 2$, $(a)$ implies that
$$
\aligned
\sup_{x\in s^G (x_\circ)}&\modo{F(x)}=2\sup_{r_1,\,r_2,\,r_3\in\Bbb Z}
\modo{(\alpha^d-1)(r_1 a_1 + r_2 a_2 + r_3 a_3)
+ i(\alpha^d+1)}^{-1}\\
&=2\left [\inf_{r_1,\,r_2,\,r_3\in\Bbb Z}
\modo{(\alpha^d-1)(r_1 a_1 + r_2 a_2 + r_3 a_3)
+ i(\alpha^d+1)}\right ]^{-1}
\le {2\over \alpha^d+1} < {2\over 3}\,,
\endaligned
$$
which proves $(b)$.
\hfill $\square$
\enddemo
\pagebreak 

The $H^\infty(X)$-hull $\widehat Y$ of a subset $Y\subseteq X$
may be found as follows:
$\widehat Y = \widehat {{\pr}_{\Bbb H} Y}\times\Bbb C$,
where ${\pr}_{\Bbb H}\colon \,X=\Bbb H\times\Bbb C\to\Bbb H$ is the natural
projection and $\widehat{{\pr}_{\Bbb H} Y}$ is the
$H^\infty(\Bbb H)$-hull of the subset ${\pr}_{\Bbb H} Y\subseteq\Bbb H$.
In view of Proposition 4.4$(b)$, we would like to pose
the following question.
\medskip

\remark{\bf 4.5. Question} For which subsets
$\Gamma\subseteq G\setminus\{\bold 1\}$
\smallskip

\noindent $(*)$ {\sl any point $x\in X$ is contained in the $H^\infty(X)$-hull
$\widehat{\Gamma (x)}$ of its $\Gamma$-orbit $\Gamma (x)$}?
\endremark
\medskip

Proposition 4.6 below provides examples of
subsets $\Gamma\subseteq G\setminus\{\bold 1\}$
with the property $(*)$.
\medskip

Recall that the elements $g_1\,,g_2\,,g_3$ form a free basis of the normal
subgroup $G_0\cong\Bbb Z^3$ in $G$.
Any subgroup $H\subseteq G_0$ is a
free abelian group of rank $\rk\,H\le 3$.

\proclaim{4.6. Proposition} Let $H\subseteq G_0$ be a subgroup
of rank $\rk\,H\ge 2$, and let $\Gamma\subseteq H$
be the complement of a finite subset\footnote
{The statement of the proposition is trivial if
$\bold 1\in\Gamma$;
however, it is meaningful whenever $\bold 1\in S$.}
$S\subset H$.
Then $x\in\widehat{\Gamma(x)}$ for any $x\in X$. In particular,
$x\in \widehat{G(x)\setminus\{x\}}$ for any $x\in X$.
\endproclaim

\demo{Proof} By Liouville Theorem, any function
$f\in H^\infty(X)=H^\infty(\Bbb H\times\Bbb C)$ is of the form
$$
f = {\widetilde f}\circ{\pr}_{\Bbb H}, \,\,\,{\text{where}} \,\,\,
\widetilde {f}\in H^\infty(\Bbb H).
\eqnum{3}
$$
Hence, for any point $x=(z,w)\in \Bbb H\times\Bbb C$
and any element $h=g_1^{r_1} g_2^{r_2} g_3^{r_3}\in H$,
we have
$$
hx=(z + r_1 a_1 + r_2 a_2 + r_3 a_3\, , w + r_1 b_1 + r_2 b_2 + r_3 b_3)
\eqnum{4}
$$
and
$$
f(hx)=\widetilde{f} (z + r_1 a_1 + r_2 a_2 + r_3 a_3)\,.
\eqnum{5}
$$
When $h$ runs over $H$ (resp. over the complement $\Gamma=H\setminus S$),
the corresponding vector
${\bold  r}=(r_1,\, r_2,\, r_3)\in\Bbb Z^3$ in (4) and (5)
runs over a sublattice $\widetilde H\subseteq\Bbb Z^3$ isomorphic to $H$
(resp. over the complement $\widetilde\Gamma
=\widetilde H\setminus\widetilde S$
of a finite subset $\widetilde S\subset\widetilde H$); in particular,
$\rk\,\widetilde H = \rk\,H\ge 2$. Since $\widetilde f$ is
a continuous function, it follows from (3), (5) and Lemma 4.2$(c)$ that
$$
f(x) = \widetilde{f}(z)\in\overline{f(\Gamma(x))}
$$
(the closure in $\Bbb C$). Therefore,
$\modo{f(x)}\le \sup_{y\in\Gamma(x)} \modo{f(y)}$,
and hence $x\in\widehat{\Gamma (x)}$.
\hfill $\square$
\enddemo

\remark{\bf 4.7. Remark} Despite Proposition 4.4$(b)$, the following fact
holds\footnote{cf. Proposition 1.10$(a)$.}:
\smallskip

\noindent {\sl For any integer $d\ne 0$ and for any
$x=(z,\,w)\in X=\Bbb H\times\Bbb C$,
the $s^G$-orbit $s^G(x)$ is a uniqueness set for bounded holomorphic
functions on $X$.}
\smallskip

That is, $f=0$ whenever $f\in H^\infty(X)$ and $f\mid s^G(x)=0$. Indeed,
any $f\in H^\infty(X)$ is of the form (3); hence, $f\mid s^G(x)=0$ implies
$\widetilde{f}\mid \pr_{\Bbb H}\,\left[s^G(x)\right] = 0$.
However, it follows from Lemmas 4.1 and 4.2(c) that the point
$\alpha^d z\in\Bbb H$ is a limit point of the set
${\pr}_{\Bbb H}\,\left[s^G(x)\right]=\left\{\alpha^d z
+ (\alpha^d-1)(r_1 a_1 + r_2 a_2 + r_3 a_3)
\mid \,r_1,\,r_2,\,r_3\in\Bbb Z\right\}
\subset\Bbb H$.
Thus, $\widetilde {f}=0$, and so $f=0$.
\endremark

\enddocument
\end